\documentstyle[12pt,twoside]{article}
\def\Ann{\mbox{Ann\,}}

\def\deg{\mbox{deg\,}}
\def\det{\mbox{det\,}}
\def\dim{\mbox{dim\,}}
\def\ext{\mbox{Ext\,}}
\def\gd{\mbox{gl. dim\,}}

\def\id{\mbox{id\,}}
\def\Im{\mbox{Im\,}}
\def\ker{\mbox{ker\,}}

\def\lm{\mbox{lm\,}}

\def\rate{\mbox{rate\,}}

\def\tor{\mbox{Tor\,}}

\def\d{\displaystyle}
\def\o{\overline}

\newtheorem{definitia}{Definition}
\newtheorem{prop}{Proposition}[section]
\newtheorem{theorema}[prop]{Theorem}

\newtheorem{corollarium}[prop]{Corolary}
\newtheorem{quest}[prop]{Question}

\def\proof{\noindent {\bf{Proof} }}

\begin{document}

\author{
        Dmitri Piontkovski\\
\footnotesize  Central Institute of Economics and Mathematics\\
\footnotesize Nakhimovsky prosp. 47, Moscow 117418,  Russia\\
\footnotesize e-mail:\
piont@mccme.ru}

%\date{}

\title{Noncommutative Koszul filtrations%
%\thanks{ Particially supported by Russian foundation of basis research project
%02-01-00468}
}

\maketitle

%{\centerline{\it Preliminary version. Comments are appreciated.}}

%\noindent

%\bigskip

%We study associative graded algebras which have a ``complete flag''
%of cyclic modules with linear free resolutions, i.e.,
%algebras over which
%there is a cyclic Koszul module with every admissible number of relations
%(from zero up
%to the number of generators of the algebra).
%Commutative algebras with the same property has been studied in several
%papers~\cite{crv, ctv,i-kos, conca,con2}.
%Here we present a non-commutative version.

\section*{Introduction}

A standard associative graded algebra $R$ over a field $k$ is called Koszul
if $k$ admits a linear resolution as an $R$-module. A (right) $R$-module $M$
is called Koszul if it admits a linear resolution too.

Here we study a special class of Koszul algebras --- roughly say,
algebras having a lot of Koszul cyclic modules. Commutative
algebras with similar properties (so-called {\it algebras with Koszul
filtrations}) has been studied in several
papers~\cite{crv, ctv,i-kos, conca,con2}.

A cyclic right $R$-module $M = R/J$ is Koszul if and only if its
defining ideal $J$ is Koszul and generated by linear forms.
So, we may deal with degree-one generated Koszul ideals
%generated by linear forms
instead of cyclic Koszul modules.
A chain $0=I_0 \subset I_1 \subset \dots \subset I_n = \overline R$
of right-sided degree-one generated ideals in an algebra $A$ is called
{\it Koszul flag} if every ideal $I_j$ is a Koszul module.
Every algebra $R$ with Koszul flag is Koszul, and the Koszul algebras
of the most important type, PBW-algebras, always contain Koszul flags
(Theorem~\ref{flag-PBW}).

In the case of commutative algebra $R$, a natural way
to find a Koszul flag in the algebra (and so, to prove its Koszulness)
is called {\it Koszul filtration}.
This concept is introduced and studied in
several papers of Conca and others~\cite{crv, ctv,conca}.
They found that the most of quadratic algebras
occuring in algebraic geometry (such as coordinate rings of canonical embeddings
of general projective curves and of projective embeddings of abelian
varieties) admits Koszul filtrations.
The main purpose of this paper is to extend their theory
to non-commutative algebras.

Namely, a set ${\bf F}$ of degree-one generated right-sided ideals in $R$
is called Koszul filtrations  if $0 \in {\bf F}, \overline R \in {\bf F}$,
and
for every $0 \ne I \in {\bf F}$
    there are $I \ne J \in {\bf F}, x \in R_1$
    such that $I = J + x R$ and the ideal $(x:J) :=
     \{ a \in R | xa \in J \}$ lies
    in ${\bf F}$.
As far as in the commutative case, every algebra admitting Koszul
filtration is Koszul (Proposition~\ref{koz.fil.}).
Every Koszul filtration containes a Koszul flag
(Corollary~\ref{fil con flag}).

For example, every monomial quadratic algebra has a Koszul filtration
(see the end of section~\ref{sect_filt} below). Another example
is an algebra with  generic quadratic relations (if either the number
of relations or $\dim R_2$ is less than $\dim R_1$, see
Theorem~\ref{generic}).

We also describe non-commutative algbras with smollest possible
Koszul filtration coinsiding with a Koszul flag,
so-called {\it initially Koszul}.
In the commutative case, algebras with the same properties
must have finite Groebner basis of commutative
relations~\cite{i-kos,ctv, crv}. We prove that,
in general, such algebras are also PBW, and we describe
their Groebner bases in Theorem~\ref{init}.
As a corollary we obtain that any algebra
with single quadratic relation (in particular, regular algebra
of global dimension two) is initially Koszul (corollary~\ref{one rel}).
Also, it follows that all tensor and semi-tensor
products of initially Koszul algebras are initially Koszul
too (Corollary~\ref{semi-tenzor}).

The relations between our classes of Koszul algebras
looks as follows
%given by the following diagram
(where arrows denote inclusions):
$$
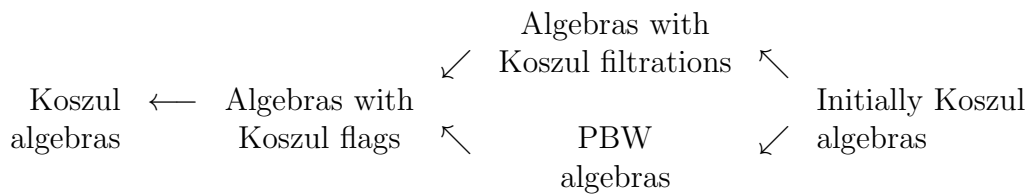
\begin{array}{rcccccl}
    & & &             & \mbox{Algebras with}       & & \\
    & & & \swarrow    & \mbox{Koszul filtrations} & \nwarrow & \\
 \mbox{Koszul}
    & \longleftarrow
	&  \mbox{Algebras with}
		    &   & & & \mbox{Initially Koszul} \\
 \mbox{algebras}
    & &  \mbox{Koszul flags}
		    & \nwarrow & \mbox{ PBW } & \swarrow
						 & \mbox{algebras} \\
   & & & & \mbox{ algebras } && \\
\end{array}
$$

Notice that the class of PBW algebras does not contain
the class of algebras with Koszul filtartions (even if we consider
commutative algebras only): this follows from~\cite[Section~3]{crv}).

\medskip

Also, we consider a generalization of the concept of Koszul filtration
for non-quadratic algebras so-called {\it generalized Koszul filtartion},
or {\it rate filtration} (see Definition~\ref{rate_def}).
In the commutative case, it is introduced
in~\cite{cnr}.  Its definition is close to the definition of Koszul 
filtration, but the ideals are not assumed to be degree-one generated.
Namely, for an ideal $I \triangleleft R$, let $m(I)$ denotes 
the maximal degree of the generators of $I$.
Then a set ${\bf F}$ of finitely generated right-sided ideals in $R$
is called rate (or generalized Koszul) filtration if 
$0 \in {\bf F}, \overline R \in {\bf F}$,
for every $0 \ne I \in {\bf F}$
    there are ideal $I \ne J \in {\bf F}$ and a homogeneous element
    $x \in I$ such that $I = J + x R, m(J) \le m(I)$, and the ideal
    $(x:J) =  \{ a \in R | xa \in J \}$ lies in ${\bf F}$.

The rate filtration is called of degree $d$, if $m(I) \le d$
for every $I \in {\bf F}$

Every algebra $R$ with rate filtration of finite degree has finite
rate (\cite{cnr}, Proposition~\ref{fil-rate}),
that is, $H_{ij}(R) = 0$  for all $j > di$ for some $d > 0$: this mens
that its homological properties are close (in a sense) to properties
of Koszul algebras. Also, we consider a class of examples
(so-called {\it algebras of restricted processing}~\cite{pi2}),
namely, a class of coherent algebras with finite Groebner basis, 
which includes all finitely presented monomial algebras.

In general, if an algebra contains a large rate filtration, 
its properties are close to the properties of coherent algebras.
In this terms, coherence means that all finitely generated 
right-sided ideals in $R$ form a rate filtration 
(Proposition~\ref{rate-coh}).

\medskip

One of misterious problems in Koszul algebras theory is to describe
the Hilbert series of Koszul algebras.
For a Koszul algebra $R$, its Yoneda algebra $\ext_R(k_R,k_R)$
coinside with {\it dual} quadratic algebra
$R^! = \bigoplus_{i} \ext^{i,i}_A(k,k)$.
The Euler characteristics for the minimal
free resolution of the trivial module $k_R$ leads to the following
relation of Hilbert series
({\it Fr\"oberg relation}):
\begin{equation}
\label{frob}
			R(z) R^!(-z) = 1.
\end{equation}
It follows that
the formal power series $R(-z)^{-1}$
has no negative coefficients.

For a long time, the following question remained opened~\cite{bac}:
 does relation~(\ref{frob}) imply the Koszul property of a given
quadratic algebra $R$?
The negative answer were obtained in 1995 independently by
L.~E.~Positselski~\cite{pos} and J.--E.~Roos~\cite{roos}. Moreover,
among all quadratic algebras with fixed Hilbert series
$R(z)$ and $R^!(z)$ there may be both Koszul and non-Koszul
algebras~\cite{pion}.

One of the most effective result in this direction is due to
Polishchuk and Positselski~\cite{pp}: thay discovered that
 for a given integer $n$, the set of
all Hilbert series of $n$--generated Koszul algebras is finite.
The same is true (for any fixed pair of integers $(p,q)$)
for the set of Hilbert series of Koszul modules
over such algebras with $p$ generators of degree $q$.

Also, they stated a natural conjecture that Hilbert series
of Koszul algebras are rational functions.
By~\cite{gov}, it is true for monomial (and so, for PBW)
algebras. By~\cite{bac, uf1} (see also~\cite{ufn}), the same holds for
quadratic algebras with at most two relations. It is proved in~\cite{dav}
that $R$-matrix Hecke algebras have rational Hilbert series too
(the Koszulness of these algebras is also proved in~\cite{wam}).
Here we prove that all algebras with Koszul filtrations
as well as all algebras with finite rate filtrations
have rational Hilbert series (Theorem~\ref{rate ratio}).

The paper is organized as follows.
In section~\ref{sect1}, we introduce our notation and give the definition of
Koszul flags. We prove that any PBW algebra containes a Koszul flag
and formulate a property of such flags (Proposition~\ref{flag property}).
In section~\ref{sect_filt}, we give the definition of Koszul filtration
and prove that any ideal which is a member of   Koszul filtration is
a Koszul module. Also, we show that any monomial quadratic algebra
has a Koszul filtration.
In section~\ref{sectgen}, we find Koszul flags and Koszul filtrations
in quadratic algebras with generic relations.
In the next section,
%~\ref{sectinitially},
we define and study initially Koszul algebras.
In section~\ref{rate}, we consider the algebras with rate filtration as
a generalization of the concept of algebras with Koszul filtrations.
In the next section, we prove that
any algebra with Koszul filtration or with finite rate filtration
has rational Hilbert series,
 and give  estimations for degrees of numerator and denominator.
Finally, in section~\ref{sectquest} we give a list of some open problems
of this theory.

\medskip

I am grateful to Leonid Positselski for helpful discussions.

\section{Koszul flags}

\label{sect1}

%Let  $k$ be a field.
We will call a vector $k$--space, a $k$--algebra, or
 $k$--algebra  module {\it graded,} if it is ${\bf Z}_+$--graded
and finite--dimensional in every component.
%For a such space  $V$ (in particular, $V$ may be an algebra or a module)
%we denote by $V(x)$ its Hilbert series
%$\sum\limits_{i \ge 0} \dim V_i x^i$.
A graded algebra $R= R_0 \oplus R_1 \oplus \dots$
is called {\it connected,} if its zero component $R_0$
is $k$; a connected algebra is called {\it standard,}
if it is generated by $R_1$ and a unit.
All algebras below are assumed to be standard, all modules and vector spaces
are assumed to be graded. All modules are right-sided.
%All inequalities between Hilbert series are coefficient--wise, i.~e., we write
%$\sum_i a_i t^i \ge \sum_i b_i t^i$ iff $a_i \ge b_i$ for all $i \ge 0$.

For a module $M$ over a $k$--algebra $R$, we will denote by $H_i (M)$
the graded vector space $\tor_i^R (M,k)$. By $H_i A$ we will denote
the graded vector space $\tor_i^R (k,k)$.

\begin{definitia}
A module $M$ is called

--- linear (of degree $d$), if it is generated in degree $d$, i.e.,
    $H_0(M)_j = 0$ for $j \ne d$;

--- quadratic, if it is lenear of degree $d$ and all
    its relations may be choosen in degree $d+1$, i.e.,
    $H_1(M)_j = 0$ for $j \ne d+1$;

--- Koszul, if it has linear free resolution, i.e.,
    $H_i(M)_j = 0$ for all $i\ge 0, j \ne i+d$.
\end{definitia}

Notice that, unlike the original definition in~\cite{bgs},
we do not assume that
quadratic or Koszul module is generated in degree zero.

\begin{definitia}{\cite{pri}}
An algebra $R$ is called Koszul if the trivial module $k_R$ is Koszul, i.e.,
every homology module $H_i R$  is concentrated in degree $i$.
\end{definitia}

In this paper
 we consider Koszul algebras having a lot of cyclic Koszul modules.

%Namely,
Let $\Phi = \{ I_0, \dots, I_n \}$ be a family of right-sided ideals of
$n$-generated stangard algebra $R$ such that  $I_0 = 0$,
$I_n = \overline R$ (maximal homogeneous ideal), and for all $0< k \le n$
$I_{k+1} = I_k + x_k,$ where $x_k \in R_1.$ We will also denote such flag
of degree one generated ideals
%this family
by $(x_1, \dots, x_n)$ and call it {\it linear flag}.
For any linear flag $\Phi = (x_1, \dots, x_n)$, the algebra $R$ is
generated by the linear forms $x_1, \dots, x_n$, and the flag $\Phi$
itself is uniqally defined by the complete flag of subspaces of $R_1$.

\begin{definitia}
A linear flag $\Phi$
is called a {\it Koszul flag} if,
 for every $0< k \le n$, the ideal $I_k$ is a Koszul module.
%    $N_k = \{ a\in R | x_k a \in I_{k-1} \}$ is generated by linear forms.
\end{definitia}

Let $S = \{x_1, \dots, x_n\}$ be a minimal system of generators of $R_1$.
Then $R \simeq F/I$, where $F = T(R_1)$ is a free associative algebra and
$I$ is a homogeneous two-sided ideal in $F$. Let us fix a linear order "$<$"
on $S$. The algebra $R$ is called {\it Poincar\'e--Birkhoff--Witt}
({\it PBW} for short) w.~r.~t. the order "$<$"
if the ideal $I$ has quadratic Groebner basis
with respect to the degree--lexicographical order  derived from $"<"$.

It is well known~\cite{pri} that any PBW algebra is Koszul.
Notice that every commutative PBW algebra has a quadratic Groebner basis
of {\it commutative }  relations; however, there are examples of
non-PBW commutative
algebras with quadratic commutative Groebner bases
%which are not PBW
(for discussions on this subject, see, e.g.,~\cite{pp}).

\begin{theorema}
\label{main}
\label{flag-PBW}
Let $R$ be a standard algebra which is minimally
generated by linear forms $x_1, \dots, x_n $.
If $R$ is a PBW algebra w.~r.~t. the order
   $x_1 < \dots < x_n$ on generators, then the
family $(x_1, \dots, x_n)$ forms a Koszul flag.
\end{theorema}

%{\bf Proof of Theorem~\ref{flag-PBW}.}
\proof

We shall prove that any ideal $I_k$ is a Koszul module.

Let $R = S/I$ as before, let $G = \{ g_1, \dots, g_r \}$ is
the reduced Groebner basis of $I$, and let $m_1, \dots, m_r$ be the leading
monomials of its members. The minimal resolution of the trivial module $k$
may be given by a construction of Anick~\cite{an2} (see also~\cite{ufn}).
In this construction, the vector space $H_i (R)$ is identified with the span
of all {\it chains} of degree $i$, i.~e., the monomials in $F$ whose
degree two submonomials lie in the set $\{ m_1, \dots, m_r\}$.
Denoting by $Ch_i$ the linear
span of such chains, we have that the minimal free resolution of $k$
has the form
\begin{equation}
\label{an_res}
0 \longleftarrow k
  \longleftarrow A
  \stackrel{d_1}{\longleftarrow} Ch_1 \otimes A
  \stackrel{d_2}{\longleftarrow} Ch_2 \otimes A
  \longleftarrow \dots
\end{equation}

For $k\le n$, let $Ch_i^k$ be the span of chains of degree $i$, whose first
letter is one of $x_1, \dots, x_k$. It follows from the
construction~\cite{an2}
that $d_i (Ch_i^k \otimes A) \subset Ch_{i-1}^k \otimes A$.
This means that there is a subcomplex ${\cal F}^k$
of the resolution~(\ref{an_res}):
$$
  {\cal F}^k:  0 \longleftarrow x_1 A + \dots + x_k A
      \stackrel{d_1^k}{\longleftarrow} Ch_1^k \otimes A
      \stackrel{d_2^k}{\longleftarrow} Ch_2^k \otimes A
  \longleftarrow \dots,
$$
where $d_i^k = d_i |_{Ch_i^k}$. In the other words, there is a filtration
${\cal F}^1 \subset \dots \subset {\cal F}^n$ on the complex~(\ref{an_res}).

It is not hard to see that the splitting homomorphisms of vector
spaces $i_j : \ker d_{j-1} \to Ch_{j} \otimes A$~\cite{an2}
are filtered too.
Since $d_j i_j = \id |_{\ker d_{j-1} }$, we obtain
$$
d_j^k ( i_j |_{\ker d_{j-1}^k} ) = (d_j i_j) |_{\ker d_{j-1}^k }
      =  \id |_{\ker d_{j-1}^k }.
$$
Thus the complex ${\cal F}^k$ is acyclic and therefore gives
a linear free resolution of the ideal $I_k = x_1 A + \dots + x_k A$.
This means that the module $I_k$ is Koszul.

\bigskip

\begin{prop}
\label{flag property}
If a family $(x_1, \dots, x_n)$ is a Koszul flag in
an algebra $R$, then (in the above notations)
all the ideals $N_k = \{ a\in R | x_k a \in I_{k-1} \}$
are generated by linear forms, and all the modules
 $N_k$ and $I_k/I_{k-1}$ are Koszul.
\end{prop}

\proof

Consider a module $M_k = I_k / I_{k-1} $. Since it is generated at degree
one, we have $H_i (M_k)_j = 0$ for $j \le i$.
The exact sequence
$$
  0 \to I_{k-1} \to I_k \to M_k \to 0
$$
leads to the long exact sequence
\begin{equation}
\label{M_k}
   \dots \leftarrow H_{i-1} (I_{k-1})_j
         \leftarrow H_i (M_k)_j
	 \leftarrow H_i (I_k)_j
	 \leftarrow  \dots
\end{equation}
For $j>i+1$, the right and left terms are zero
(since $I_{k-1}$ and $I_k$ are Koszul).
Hence $H_i (M_k)_j = 0$
for $j \ne i+1$, that is, $M_k$ is Koszul.

Now
%, to prove that the statements about the ideals $N_k$,
it remains to note that
$M_k \simeq R/N_k[-1]$, so, $N_k$ is generated by linear forms and Koszul.

%By the same way, the module $N_k$ is concentrated at degrees greater or
%equal to one, so $H_i (N_k)_j = 0$ for $j < i+1$.
%From the short exact sequence
%$$    0 \to \Omega I_{k-1} \to \Omega I_k \to N_k[-1] \to 0
%$$
%(see the proof of Proposition~\ref{koz.fil.}), we obtain
%\begin{equation}
%\label{N_k}
%  0 \leftarrow  H_0(N_k)_{j-1} \leftarrow H_1(I_k)_{j}
%     \leftarrow H_1(I_{k-1})_{j}
%\end{equation}
%and
%$$
%    H_{i} (I_{k-1})_j
%         \longleftarrow H_i (N_k)_{j-1}
%	 \longleftarrow H_{i+1} (I_k)_j
%$$
%for $i \ge 1$. Since right and left terms are zero for $j >i$,
%we have $H_i(N_k)_j = 0$ for $j \ne i+1$, so $N_k$ is Koszul too.

\section {Noncommutative Koszul filtrations}

\label{sect_filt}

%In this section, all our results and methods are essentially
%the same as in~\cite{ctv, crv} where the only case of commutative algebras
%has been considered.

A natural way leading to Koszul flags is given by the following
concept~\cite{ctv, crv}.

\begin{definitia}
\label{kos.filtr.}
A family ${\bf F}$ of right-sided ideals of a standard algebra $R$ is said
to be {\it Koszul filtration} if:

1)  every ideal $I \in {\bf F}$ is generated
by linear form;

2) zero ideal and maximal homogeneous ideal $\overline R$ are in
${\bf F}$;

3) for every $0 \ne I \in {\bf F}$
    there are $I \ne J \in {\bf F}, x \in R_1$
    such that $I = J + x R$ and the ideal $N = \{ a \in R | xa \in J \}$ lies
    in ${\bf F}$.
\end{definitia}

Similarly to the commutative case~\cite{ctv, crv},  this notations is due to
the following

\begin{prop}
\label{koz.fil.}
Let ${\bf F}$ be a Koszul filtration on $R$. Then every ideal
$I \in {\bf F}$ is a Koszul module.
\end{prop}

The proof is close to the one of the commutative
version~\cite{ctv, crv}.

\proof

We have to prove that $H_i(I)_j = 0 $  for all $i\ge 0, j \ne i+1$.
We will use the induction by both $i$ and the number of generators of $I$.
If $I=0$ or $i = 0$, it is nothing to prove. Otherwise let $J,N$ be the same
as in
Definition~\ref{kos.filtr.}.

Put $x_k = x$.
Let $\{ x_1, \dots, x_{k-1} \}$ is a minimal set of linear forms generated
$J$, and let $\{ X_1, \dots, X_{k-1} \}$
(respectively,  $\{ X_1, \dots, X_k \}$   )  is the set of generators
of $H_0 (J)$  (resp., $H_0 (I)$ ) such that their images in $J$ (resp., $I$)
are $x_1, \dots, x_{k-1}$ ($x_1, \dots, x_{k}$).
Consider the first terms of minimal free resolutions of $J$ and $I$:
$$
   0 \longleftarrow J \longleftarrow H_0 (J) \otimes R \longleftarrow
	   \Omega (J) \longleftarrow 0,
$$
$$
   0 \longleftarrow I \longleftarrow H_0 (I) \otimes R \longleftarrow
	   \Omega (I) \longleftarrow 0,
$$
where the syzygy modules are
$\Omega (J) = \{ X_1 \otimes a_1 + \dots + X_{k-1} \otimes a_{k-1} |
		 x_1 a_1 + \dots + x_{k-1} a_{k-1} = 0
              \}$
and
$\Omega (I) = \{ X_1 \otimes a_1 + \dots + X_{k} \otimes a_{k} |
		 x_1 a_1 + \dots + x_{k} a_{k} = 0
              \}$.
We have $\Omega (J) \subset \Omega (I) $ and
$ \Omega (I) / \Omega (J) \simeq \{ X_k \otimes a | x_k a \in J\}
  \simeq N [-1].
$
Tensoring the short exact sequence
$$
  0 \longrightarrow \Omega (J) \longrightarrow
     \Omega (I)  \longrightarrow N [-1] \longrightarrow 0
$$
by $k$, we obtain from the long exact sequence of $\tor$'s that
the sequence
$$
  H_i (J)_j \to H_{i} (I)_j \to H_{i-1} (N)_{j-1}
$$
is exact for all $i \ge 1$.
By induction, the right and left terms vanish for $j \ne i+1$;
so the middle term vanishes too.

\medskip

\begin{corollarium}
\label{fil con flag}
Every Koszul filtration containes a Koszul flag.
\end{corollarium}

{\noindent{\bf Example: monomial algebras}}

Let $M$ be an algebra with monomial quadratic relations.
Let $X = \{ x_1, \dots, x_n\}$ be the set of its generators.
Consider the set ${\bf F}$ of all right-sided ideals in $M$
generated by subsets of $X$. We claim that ${\bf F}$ forms a
Koszul filtration. (Algebras with such a property are called
{\it strongly Koszul}. Another example is a commutative
monomial algebra~\cite{HHR}.)

Indeed, let $I = J + x_j M$ be two ideals in ${\bf F}$. We have to prove
that the ideal $N = \{ a \in M | \, x_j a \in J \}$ lies
    in ${\bf F}$.
Consider any linear flag
 $\Phi = (x_{\sigma 1}, \dots, x_{\sigma n})$ (where $\sigma$
is a permutation) such that $I$ and $J$ are its elements.
By Theorem~\ref{flag-PBW}, it is a Koszul flag.
It follows from Proposition~\ref{flag property} that the ideal $N$
is generated by linear forms. Since the algebra $M$ is monomial,
$N$ is generated by a subset of $X$, so $N \in {\bf F}$.

\section {Quadratic algebras with generic relations}

\label{sectgen}

Let $R$ be a quadratic algebra with $n$ generators and $r$
generic quadratic relations. It is well-known~\cite{an1} that
$R$ is Koszul iff either $r \le n^2/4$ or $r \ge 3n^2/4$.
Let $x_1, \dots, x_n$ be a generic set of generators of $R$.
Consider the generic linear flag $\Phi = (x_1, \dots, x_n)$.

\begin{theorema}
\label{generic}
If either $r < n$ or $r > n^2 - n$, then the generic linear flag $\Phi$
is a Koszul flag and is a subset of some Koszul filtration ${\bf F}$.
If $r<n$, one can take the filtration  ${\bf F}$ to be finite.

If $n \le r \le n^2 - n$, then  $\Phi$ is not a Koszul flag, and so, is not
a part of a Koszul filtration.
\end{theorema}

\proof

Let $R = F/I$, where $F$ is a free algebra generated by  $x_1, \dots, x_n$
and $I$ be a two-sided ideal generated by $r$ generic quadratic forms
$f_1, \dots, f_r$. This means that
 $f_j = \sum\limits_{i \le n} x_i l_j^i$,
where $l_j^i$ are generic linear forms.

First, consider the case $r<n$. Then $R$ has global dimension two, that is,
minimal free resolution of the trivial module $k_R$
has the form
  \begin{equation}
  \label{res}
	0 \longleftarrow k \stackrel{d_0}{\longleftarrow}
        R \stackrel{d_1}{\longleftarrow}
        H_1 \otimes R \stackrel{d_2}{\longleftarrow}
        H_2 \otimes R \longleftarrow 0,
  \end{equation}
where the vector space $H_1$ is the span of indeterminates
$\tilde x_1, \dots, \tilde x_n$,
and the vector space $H_2$ is the span of indeterminates
$\tilde f_1, \dots, \tilde f_r$. Differentials here send
$\tilde x_i$ to $x_i$ and $\tilde f_i$   to $f_i$.

Denote the ideals of flag $\Phi$ by $I_t = \{ x_1, \dots, x_t\} R$
and $I_0 = 0$.
We claim that the ideals $N_t = \{ a\in R | x_t A \in I_{t-1} \}$
vanish for all $1 \le t < n$.

Indeed, assume the converse. Then $N_t \ne 0$ for some $t < n$.
By definition, it means
that there are elements $c_1, \dots, c_t \in R$, not all zeroes, such that
for $t$ generic elements $x_1, \dots, x_t$ we have
$$
\sum\limits_{i \le t} x_i c_i = 0.
$$ We may assume that
all the elements $c_i$ have the minimal possible degree
(say, $d$) of all such $t$-ples for all $t < n$.

Let us denote $f = \sum\limits_{i \le t} \tilde x_i c_i \in H_1 \otimes R$.
In the resolution~(\ref{res}),
we have $d_1 (f) =
\sum\limits_{i \le t} x_i c_i = 0$,
so, $f$ lies in the module $\Im d_2$.

 It follows that, for some $a^1, \dots, a^r \in R$,
we have
$$
f = \sum\limits_{j \le r} d_2(\tilde f_j) a^j
 = \sum_{j \le r, i\le n} \tilde x_i l_j^i a^j.
$$
Since $\deg f = d+1$, for every $j$ we may assume that
$\deg a^j  = d -1$.
Taking the projection on $\tilde x_n \otimes R $,
we have
$$
 \sum_{j \le r}  l_j^n a^j = 0.
$$
We obtain the linear relation between $r < n$ generic elements
$l_1^n, \dots, l_r^n$ with coefficients of degree $d-1$. This contradicts to
the minimality of $d$.

Now, consider the ideal $N_n$. Proceeding as above, we have that an element
$a = a_n$ lies in $N_n$ if and only if there are elements
$a_1, \dots, a_{n-1}$ such that $\sum\limits_{i \le n} x_i a_i = 0$.
This means that $N_n$ is isomorhical to the projection of the module
$\Im d_2$ on the component $\tilde x_n \otimes R $. So, $N_n$ is generated
by $r$ generic linear forms $l_1^n, \dots, l_r^n$.

It has been proved that, for any
$t < n$ generic linear forms $x_1, \dots, x_t$,
all the ideals $N_i = \{ a | x_i a \in (x_1, \dots, x_{i-1})R \}$ vanish.
So, the same is true for the linear forms $l_1^n, \dots, l_r^n$.
Thus, the $n+r +1$ ideals
$$
J_0 = 0 = I_0, I_1, \dots, I_n, J_1 = l_1^n R, \dots,
          J_r = l_1^n R + \dots + l_r^n R
$$
form a Koszul filtration. For completeness, we show how the filtration
looks in terms
of Definition~\ref{kos.filtr.}:
$$
\begin{array}{c|c|c}
I &                 J &                N \\
I_t, 1 \le t< n     &   I_{t-1}     &     0 \\
I_n = R_+           &   I_{n-1}     &   J_r \\
J_t, 1 \le  t \le r &  J_{t-1}      &   0  \\
\end{array}
$$

Now, consider the case $r > n^2-n$. We will construct
the filtration ${\bf F}$ starting from the flag $\Phi$.

We have $s := \dim R_2 = n^2 - r < n$
and $R_3 = 0$. Since $R_2 = x R_1$ for all generic $x \in R_1$,
for $2 \le t \le n$ we have
$N_t = \{ a | x_t a \in I_{t-1} \} \supset \{ a | x_t a \in x_1 R \} = R_+$.
Consider the ideal $N_1 = \Ann_R x_1$.
It containes
(and so, generated by) $n-s$ generic linear forms, say,
$x^1_1, \dots, x^1_{n-s}$. Let us add the ideals
$I_t^1 = \sum\limits_{i \le t} x^1_i R, t \le n-s$, to the filtration.
To make the filtration Koszul,  we need to add also the ideals
$N^1_t = \{ a | x^1_t a \in I^1_{t-1} \}$. By the above reasons, we have
$N^1_t = 0$ for $t \ge 2$, but the ideal $N_1^1 = \Ann_R x_1^1$
is generated by
another set of $n-s$ generic linear forms, say, $x^2_1, \dots, x^2_{n-s}$.
Proceeding as well, we get the infinite Koszul filtration ${\bf F}$.

\medskip

Now, let us prove the second (negative) part of the Theorem.
Since the albegra $R$ is not Koszul for $n^2/4 < r < 3n^2/4$,
we have to consider two cases: $n \le r \le n^2/4$ (then $\gd R =2$)
and $3n^2/4 \le r \le n^2 -n$ (then $R_3 = 0$).

Let $r = qn +p$, where $0 \le p < n$.
Taking Gaussian elemination of the monomials in the quadratic forms
$f_1, \dots, f_r$,
we may assume that every monomial $x_i x_j$ (where lexicographically
$(i,j) \le (q+1,p) $) appears only at $f_{n(i-1) + j}$.
This means that, under the lexicographical order with $x_1 > \dots > x_n$,
we have $\lm f_{n(i-1) + j} = x_i x_j$.

Assume that $\Phi$ is a Koszul flag.

Consider the Artinian case:  $3n^2/4 \le r \le n^2 -n$.
By Proposition~\ref{flag property}, all the modules $N_t$ are generated
by linear forms. Since $R_2 \subset N_1$, we have that there is at least one
nonzero
linear form (say, $y$) at $N_1 = \Ann_R x_1$.

If $y = \sum_{i\le n} \lambda_i x_i$, we have that in the free algebra $F$
$$
    x_1 \sum_{i\le n} \lambda_i x_i = \sum_{i\le n} \lambda_i x_1 x_i
    \in k \{ f_1, \dots, f_r\}.
$$

Taking Gaussian elimination of monomials, we may assume that
$f_i = x_1x_i + g_i, i \le n$,
where the monomials
$x_1x_i, 1\le i \le n$ do not appear at the quadratic
forms $f_{n+1}, \dots, f_{r} $ and $g_1, \dots, g_n$.
This get the equality of elements of $F$
$$
    \sum_{i\le n} \lambda_i x_1 x_i = \sum _{i\le n} \lambda_i f_i +
     \langle \mbox{linear combination of other $f_i$'s}  \rangle.
$$
Since the monomials $\lm f_i$ for $i > n$ do not appear at the left side
of the equality, they do not also appear at the right side. So,
the last addition vanishes, and we get the  equality
$
\sum_{i\le n} \lambda_i x_1 x_i = \sum _{i\le n} \lambda_i f_i
$,
or
$$
\sum_{i\le n} \lambda_i g_i = 0.
$$
This is a system of linear equations on $n$ variables $\lambda_i$
with generic coefficients. The number of equations is equal to
the number of monomials in every $g_i$, i.e., $\dim R_2 = n^2 -r \ge n$.
So, the unique solution is the zero vector. This contradicts to the choice
of $y$.

\medskip

It remains to consider the case of genric algebra of global dimesion two
with $n \le r \le n^2/4$.
By Proposition~\ref{flag property}, all the modules $N_t$ are generated
by linear forms. Taking their degree-one components, we get
$N_n = \dots = N_{n-q+1} = R_+$,
$N_{n-q-1} = \dots = N_1 = 0$, and $N_{n-q}$ is generated by $p$ generic
linear forms, so that $N_{n-q}(z) = I_p(z)$.

By the exact sequences
$$
    0 \to I_t \to I_{t+1} \to R/N_{t+1} [1] \to 0,
$$
we obtain the formulae for Hilbert series:
$$
	  I_{t+1}(z) = I_t(z) +z \left( R(z) - N_{t+1}(z) \right).
$$
Thus,
$$
  I_0(z) = 0, I_1(z) = zR(z), \dots, I_{n-q-1} (z) = (n-q-1)zR(z),
$$
$$
    I_{n-q+1} (z) = I_{n-q}(z) +z,  \dots, I_{n} (z) = I_{n-q}(z) + qz,
$$
and
$$
	    I_{n-q}(z) = I_{n-q-1} (z) +z \left( R(z) - I_p(z) \right)
		= (n-q) z R(z) - zI_p(z) .
$$
Since $I_n = R_+$,
we get
$$
   I_n(z) = R(z) - 1 = qz +(n-q) z R(z) - z I_p(z).
$$
Here $R(z) = \left( 1-nz+ rz^2 \right)^{-1} $ ("Golod--Shafarevich"),
thus
$$
    I_p(z) = z \frac{\d p+qrz}{ 1- nz+rz^2} = (zp+qrz^2) R(z).
$$

Since $\gd R =2$, the minimal free resolution of the module $R/I_p$
has the form
$$
   0 \gets R/I_p  \gets R  \gets H_0(I_p) \otimes R
          \gets H_1(I_p) \otimes R \gets 0.
$$
Here $H_0(I_p) (z) = pz$, so, the Euler characteristic leads to the formula
$$
   I_p (z) = pz R(z) - H_1 (I_p)(z) R(z).
$$
Finally, we obtain
$$
    H_1 (I_p)(z) = pz - R(z)^{-1}I_p (z) = - qrz^2.
$$
For $q >0$, the coefficient iz negative, so, it cannot be a Hilbert series.
The contradiction completes the proof.

\section {Initially Koszul algebras}

\label{sectinitially}

If a Koszul flag $(x_1, \dots, x_n)$ in an algebra
$R$ forms a Koszul filtration, it is called {\it Groebner flag}~\cite{ctv, crv};
an algebra $R$ having such a flag is said to be {\it initially Koszul}
(w.~r.~t. the sequence of generators  $x_1, \dots, x_n$)~\cite{i-kos}.

In the commutative case~\cite{i-kos,ctv, crv},
every initially Koszul algebra
has a quadratic Groebner basis w.~r.~t. the reverse lexicographical order.
In our non-commutative case, such an algebra has a quadratic Groebner basis
w.~r.~t. the ordinar degree-lexicographical order.

\begin{theorema}
\label{init}

Let $R$ be a standard algebra generated by degree one elements
$x_1, \dots, x_n$.
Then the following two conditions are equivalent:

(i) $R$ is an initially Koszul algebra with Groebner flag
$(x_1, \dots, x_n)$;

(ii) $R$ is PBW w.~r.~t. the degree-lexicographical order with
$x_1 < \dots < x_n$,   and the leading monomials of the Groebner
basis $G$  of the relations ideal has the following property:

if $x_k x_j \in \lm G$, then $x_k x_i \in \lm G$ for all $i < j$;

(iii) the monomial algebra $R' = F/\id (\lm G)$ is initially Koszul
with the same Groebner flag.

\end{theorema}

\proof

$(i) \Longrightarrow (ii) $

First, let us prove that $R$ is PBW.

Suppose that a minimal Groebner basis $G$ of $I$ includes an element
$g$ whose degree is greater or equal to 3 (where $I$ is the two-sided
ideal of relations of $R$ in the free algebra
$F = k \langle x_1, \dots, x_n \rangle$). Let
$g = x_k s_k + x_{k-1} s_{k-1} + \dots + x_1 s_1$, where $s_k \not\in I$
and $\deg s_i \ge 2$ for all $s_i \ne 0$.

Consider the ideal $N_k$. For some $r = r_k$, it is
generated by
$x_1, \dots, x_r$.
 Since $x_i \in N_k$ for $i \le r$,
we have that for every $i \le r$
there is an element %of $F$
$f^k_i = x_k (x_i + \alpha^k_{i-1}x_{i-1} +\dots + \alpha^k_{1}x_{1}) +h^k_i
\in I$, where $\alpha^k_j \in k, h^k_i \in I_{k-1}$.
Since the leading monomials $\lm {f^k_i} = x_k x_i$
of the elements $f^k_i$
%of these quadratic elements
are pairwise different, we may assume that
they are the members of our Groebner basis $G$.

For an element $a \in F$, %let us
denote by $\overline {a}$ its
image in $R$.
Since $g \in I$, it follows that  $\overline {s_k} \in N_k$, so, for some
$t_1, \dots, t_r \in F$ we have
$\o s_k = x_1 \o t_1 + \dots + x_r \o t_r$.
Also, we may assume that
%It follows that
the leading monomial $\lm s_k$ is equal to
$x_i \lm t_i$ for some $i \le r$.
Thus
$$
\lm g = x_k \lm s_k = x_k x_i \lm t_i
= \lm (f^k_i) \lm t_i.
$$
This means that  $g$ is reducible w.~r.~t. $f^k_i$, in contradiction to our
assumption that the Groebner basis $G$ is minimal. So, $R$ is PBW.

To complete the proof of the implication, we will show that
$G$ consists of the elements $f_i^k$ for all $k \le n, i \le r_k$.
Indeed, assume that there is $g \in G$ such that $\lm g = x_k x_t$
for some $t > r_k$. This means that
$ g = x_k v + x_k w + h$, where $v$ lies in the span
of $x_{r+1}, \dots, x_t$, $w$ lies in the span of
$x_{1}, \dots, x_r$, and $h$ lies in the right-sided ideal generated
by $x_1, \dots, x_{k-1}$. It follows that $\o v \in N_k$,
in contradiction to our assumption that $N_k$ is generated by
$x_1, \dots, x_r$.

\medskip

$(ii)  \Longrightarrow (i) $

It follows form~$(ii)$ that for every
$1 \le k \le n$ there is an integer $r=r_k$
such that the leading monomials of $G$ are exactly
$x_k x_i$ for all $k \le n, i \le r_k$.

 Let us show that, for every $k \le n$, the ideal $N_k$ is generated
by $x_1, \dots, x_r$. Let $\o a \in N_k$, where $a \in F$
is an irreducible element. Since
$x_k a \in$ $\{ x_1, \dots, x_{k-1} \} F + I$, it follows that
the monomial $x_k \lm a$ is reducible.
If $\lm a = x_{i_1} \dots x_{i_d}$,
this means that for some $j \le d$ there is an element $g$
of the Groebner basis $G$ such that
$\lm g = x_k x_{i_1} \dots x_{i_j}$. Thus $a \in$
$\{ x_1, \dots, x_{r} \} F + I$, so $\o a \in I_r$.

\medskip

$(i)  \Longleftrightarrow (iii) $

By definition, the sets of leading terms of Groebner basis
are coinside for the algebras $R$ and $R'$.
By the equivalence $(i)  \Longleftrightarrow (ii) $,
it follows that these algebras are initially Koszul simultaneously.

\begin{corollarium}
\label{one rel}
Suppose that the ground field $k$ is algebraically closed.
Then every quadratic algebra $R$ with one relation is initially Koszul.

For example, standard (Artin--Shelter) regular
algebras of global dimension~2~\cite{zhang}
are initially Koszul.
\end{corollarium}

\proof

Let $R = F/\id (f)$, where $F$ is a free algebra and $f$ is the relation.
Let ${\cal L}(f)$ be a set of all degree one generated right-sided ideals
in $F$ which contain $f$, and let $J$ be any minimal (by inclusion)
element of ${\cal L}(f)$.

If $J$ is a principal ideal, then $f$ has the form $xy$, where $x,y$
are linear forms. This imply that, under a linear change of generators
$x_1, \dots, x_n$
of $F$, either $f = x_n x_1$ or $f= x_1^2$. In both cases, $R$
is initially Koszul by Theorem~\ref{init}.

Now, let $J$ is generated by at least two elements.
We may assume that $x_{n-1}, x_n \in J$. It is
easy to see~cite[p.~1307]{pi3} that, up to
a linear transformation of the generators, the representation of $f$
does not contain the term $x_n^2$. That is, $f = x_n l + l' x_n +g$,
where $l,l', g$ do not depend on $x_n$.
By the minimality of $J$, $l \ne 0$. Making any linear transformation
of variables
$x_1, \dots, x_{n-1}$
sending $l$ to $x_1$, we obtain $\lm f = x_n x_1$.
By Theorem~\ref{init}, this means that $R$
is initially Koszul.

\begin{corollarium}
\label{monom dual}
Suppose that a monomial quadratic algebra $R$ is initially Koszul
with Groebner flag $(x_1, \dots, x_n)$, where $x_1, \dots, x_n$
are the monomial generators.
Then its quadratic dual algebra $R^! $ is initially Koszul with
Groebner flag $(x_n^*, \dots, x_1^*)$, where asterisque denotes
the dual element.
\end{corollarium}

\proof

In the dual basis, the dual monomial algebra has the relations
$$
\{x_i^* x_j^*  | x_i x_j \mbox{ is {\it not} a relation of }R \}.
$$
Then the Corollary follows from the
Theorema \ref{init}.

\medskip

For the next corollary, we need a generalization of the concept
of tensor product~\cite{an1}.
Let $A$ and $B$ be two standard algebras generated
by the sets of indeterminates $X$ and $Y$, and let their Groebner bases of
realations
(w.~r.~t. suitable degree--lexicographical orders) be
 $G_A$ and $G_B$, respectively. Assume that an algebra
$C$ is a quotient of the free product $A * B$.
%algebra $T (kX + kY)$.
Then $C$ is called a {\it semi--tensor product} of $A$ and $B$
iff its Groebner basis of relations $G_C$
(w.~r.~t. the  degree--lexicographical order with $X > Y$)
has the same set of leading monomials as the Groebner basis of the tensor product
$A \otimes B$, i.~e.,
$$
\lm G_C = \lm  G_A \cup \lm  G_B \cup \{ xy | x \in X, y \in Y\}.
$$

\begin{corollarium}
\label{semi-tenzor}
A semi-tensor product of initially Koszul algebras is initially Koszul.

For example, quantum polynomial rings
%$k \langle x_1, \dots, x_n | x_i x_j = q_{ij} x_jx_i\rangle$,where $q_{ij} q_{ji} =1$,
%standard skew polynomial rings, and quadratic Artin--Shelter
%regular algebras of global dimension~3 \cite{AS, ATB}
are initially Koszul.
\end{corollarium}

\section{Generalized Koszul filtrations, or Rate filtrations}

\label{rate}

\subsection{Definition and main property}

The notion of {\it generalized Koszul filtration} is
an analogue of Koszul filtration for non-quadratic commutative
algebras~\cite{cnr}. Their definition admits direct non-commutative
generalization. For briefness, we also call this filtration as
{\it rate filtration}.

For a graded right-sided ideal $I \in R$, let $m(I)$
denote the maximal degree of its homogeneous generator.

\begin{definitia}
\label{rate_fil_def}
\label{rate_def}
Let $R$ be a
standard (i.e., degree-one generated)
finitely generated graded algebra, and let ${\bf F}$
be a set of finitely generated right-sided ideals in $R$.
The family $F$ is said to be  generalized Koszul filtration, or
rate filtration,
if:

1) zero ideal and the maximal homogeneous ideal $\overline R$
belong to $F$, and

2) for every $0 \ne I \in {\bf F}$
    there are ideal $I \ne J \in {\bf F}$ and a homogeneous element
    $x \in I$ such that $I = J + x R, m(J) \le m(I)$, and the ideal
    $N = (x:J) =  \{ a \in R | xa \in J \}$ lies in ${\bf F}$.
\end{definitia}

We will say that a rate filtration
is {\bf of degree $d$} if all
its members are generated at degrees at most $d$.
Koszul filtrations are exactly rate filtrations
of degree one.

Let us recall a notations.
% of~\cite{cnr}.
For a graded finitely generated $R$--module $M_R$, put
$t_i = \max \{ j | H_i(M)_j \ne 0 \}$; if $H_i(M) = 0$, put
$t_i = 0$.
The {\it rate}~\cite{brate} of algebra $R$
is the number
$$
    \rate R = \sup_{i \ge 2} \{ \frac{t_i(k) - 1}{i-1} \}.
$$
For commutative standard algebras~\cite{abh} as well as for
non-commutative algebras with finite Groiebner basis of relations~\cite{an2}
 the rate is always finite.
Rate is equal to 1 if and only if $R$ is Koszul.
If an algebra has finite rate, then its Veronese subring of sufficiently
high order is Koszul.

The following Proposition is
originally proved for commutative
algebras~\cite{cnr}. In fact, it holds for non-commutative
ones as well. It is an analogue of Proposition~\ref{koz.fil.}
for rate filtrations instead of Koszul filtrations.

\begin{prop}
\label{fil-rate}
Let ${\bf F}$ be a rate filtration of degree $d$ in $R$.
%Suppose that  there is an integer $d$ such that
%$m (I) \le d$  for every $I \in {\bf F}$.
Then
$$
      t_i(I)  \le m(I) + di
$$
for all $i \ge 0$ and $I \in {\bf F}$;
in particular,
%it follows that
$$
      \rate R \le d.
$$
\end{prop}

\proof

Like the proof of Proposition~\ref{koz.fil.}
(and of the commutative
version from~\cite{cnr}), we proceed by induction over $i$ and on $I$
(by inclusion). First, notice that the degree $c$ of $x$ in
Definition~\ref{rate_fil_def} cannot be greater than $m(I) \le d$.
Proceeding as in Proposition~\ref{koz.fil.}, we obtain
the triple
$$
  H_i (J)_j \to H_{i} (I)_j \to H_{i-1} (N)_{j-c}.
$$
By induction, first term vanishes for $j \ge m(J)+d(i-1)$,
and the third one vanishes for $j-c \ge m(N) + d(i-2)$.
Since $m(J)\le m(I)$ and $m(N) \le d$, they both vanish for all
$j \ge m(I) + d(i-1)$, and so the middle term vanishes too.

\subsection{Examples: monomial algebras and similar constructions}

Like the commutative case~\cite{cnr}, in any standard
monomial algebra $Q$
whose relations are of degrees at most $d$, there is a
rate filtration {\bf F} of degree $d-1$: it consists of all monomial ideals
%in $Q$
generated in degrees less than $d$.

Another class of examples is the class of so-called ``algebras of
restricted processing''~\cite{pi2}, that is, algebras with finite
Groebner basis of a special kind. Let $A$ be a quotient
algebra of a free algebra $F$ by an ideal $I$ with
Groebner basis $G = \{ g_1, \dots, g_s \}$. For every element
$f \in F$, it is well-defined its normal form N(f) with respect to $G$.
For some $r \ge 0$, algebra $A$ is called {\it algebra of $r$--processing},
if for any pair $p,q \in F$ of normal monomials,
where $q = q_1 q_2, \deg q_1 \le r$,
$$
      N(p q) = N(p q_1) q_2.
$$
The simple example is a monomial algebra $A$ presented by monomials
of degree at most $r+1$.

A sufficient condition for an algebra to have this property is as follows.
Consider a graph $\Gamma$ with vertices marked by $g_1, \dots, g_s$.
An arrow $g_i \to g_j$ exists iff there is an overlap between
any {\it non-leading} term of $g_i$  and {\it leading} term of $g_j$.
If $\Gamma$ is acyclic, then  $A$ is an algebra with $r$--processing for some
$r$.

\begin{prop}
Let a standard algebra $A$  be an algebra with $r$--processing.
Then all its right-sided ideals generated in degree at most
$r$ form a rate filtration (of degree $r$).
\end{prop}

In particular, any  algebra of $1$--processing is {\it universally Koszul}
(that is, the set of all its degree-one generated
ideals forms a Koszul filtration~\cite{conca}).  Such algebras
were separately considered in~\cite{iou}.

\proof

Let $I \subset A$ is a right-sided ideal with $m(I) \le  r$
which is minimally generated by a set $X$. Let $x\in X$
be any generator of degree $m(I)$, let $X = \{x\} \bigsqcup Y$,
 and let  $J = Y A $. It is sufficient to show that
$m(N) \le r$, where $N = \{a \in A | xa \in J \}$.
This follows from~\cite{pi2}[Lemma~4].

\subsection{Rate filtrations and coherence}

Recall that an algebra $R$ is called {\it (right) coherent}
if every map $M \to N$ of two finitely generated
(right) free $R$-modules has finitely
generated kernel; other equivalent conditions (such as, every
finitely generated right-sided ideal in $R$ is finitely presented)
 may
be found in~\cite{faith,burbaki}.
Every Noetherian ring is coherent. Free associative algebra
is coherent, as well as finitely presented monomial algebras and
algebras of $R$-processing~\cite{pi2}. If the algebra $R$
is graded, it may introduce two versions of coherence,
``affine'' and ``projective'' (where all maps and modules are assumed
 to be graded): the author does not know whether these concepts are
equivalent or not.
The (projective) coherent rings may be considered as
a basic of non-commutative geometry instead of Noetherian rings~\cite{pol}.

One of equivalent definitions of coherent rings is as follows
(Chase, cited by~\cite{faith}): $R$ is coherent iff, for every finitely
generated ideal $J = JR$ and element $x \in R$, the ideal
$N = (x:J) := \{ a \in R | xa \in J \}$  is finitely generated.
This definition is close to our definition of rate filtration, as shows the
following

\begin{prop}
\label{rate-coh}

For a standard algebra $R$, the following two statement are equaivalent:

(i) $R$ is projective coherent

(ii) all finitely generated homogeneous ideals in $R$ form a rate filtration
(of infinite degree).
\end{prop}

\proof

The implication $(i) \Longrightarrow (ii) $
follows from the definition above. Let us prove
$(ii) \Longrightarrow (i)$. We will show that any finitely generated
right-sided ideal $I$ in $R$ is finitely presented.

Let us proseed by induction both by $m(I)$ and the number of generators
$n$ of $I$. Let $I$ be minimally generated by a set $X =\{x_1, \dots, x_n \} $.
We have to show that the syzygy module $\Omega = \Omega(I)$
is finitely generated, where
$$
    0 \longleftarrow I \longleftarrow
                        kX \otimes R \longleftarrow \Omega \longleftarrow 0.
$$

Let $x, J, N$ be the same as in the definition of rate filtration,
and let us suppose that $x = x_n$ and $J$ is generated by
the set $X' = \{ x_1, \dots, x_{n-1} \} $.
In the exact triple
$$
    0 \longrightarrow  \Omega(J) \longrightarrow  \Omega(I)
	\longrightarrow  \Omega(I)/ \Omega(J) \longrightarrow 0
$$
we have  $\Omega(I)/ \Omega(J) \simeq N [-1]$.
By the induction, both the first and the last modules in the triple
are of finite type, and so the middle is.

\medskip

In the same way, the following corollary on rate filtrations
offinite degree may be shown.
It generalizes the property of algebras of $r$-processing.

\begin{corollarium}
Assume that, for every sufficiently large integer $d$,
all right-sided ideals $I$ in the algebra $R$ with
$m(I) \le d$ form a rate filtration. Then
$R$ is projective coherent.
\end{corollarium}

\section {Rationality of Hilbert series}

Here we
consider a rate filtration ${\bf F}$
of finite degree $d$ such that the number of pairwise
different Hilbert series of ideals $I \in {\bf F}$ is finite.
Our main examples are as follows:

1) finite rate filtration, in particular, arbitrary  rate filtration
of degree $d$ over a finite field $k$ (where any ideal $I \in {\bf F}$
is generated by a subset of finite set $R_{\le d}$);

2) any Koszul filtration.

Indeed, it is proved in~\cite{pp} that for every Koszul algebra $R$,
for given
integer $n$ the set of Poincare series of all degree-one generated
Koszul $R$--modules   with at most $n$ generators is finite.
It follows that the set of all Hilbert series of degree-one generated ideals
is finite.

\begin{theorema}
\label{rate ratio}
Suppose that an algebra $R$ has a
rate filtration ${\bf F}$ of some degree $d$ such
that the set ${\cal H}ilb$ of all Hilbert series of ideals
$I \in {\bf F}$ is finite.
Then $R$ has rational Hilbert series, as well as every ideal $I \in {\bf F}$.

If ${\cal H}ilb$ contains $s$ nonzero elements, then the degrees
of numerators and denominators of these rational functions are not
greater than $ds$.
\end{theorema}

\proof

Let ${\cal H}ilb = \{I_0(z) = 0, I_1(z), \dots, I_s(z)\}$,
where $I_1, \dots, I_s$ are some nonzero ideals liing in ${\bf F}$.
By definition, for every  nonzero ideal $I = I_i$
there are elements $ J = J(I), N = N(I) \in {\bf F}$
such that $J \in I, J \ne I$,     and for some positive
degree $d = d(I)$ the following triple is exact:
$$
     0 \to J \to I \to R/N [-d] \to 0.
$$
Taking Euler characteristics, we have
\begin{equation}
\label{rate_hilb}
     I(z) = J(I)(z) + z^d (R(z) - N(I)(z)).
\end{equation}
Let us put $J^{(1)} := J(I)$ and $J^{(n+1)} := J\left( J^{(n)} \right)$.
Since all ideals $J^{(n)}$  are generated by subspaces of $R_{\le d}$,
the chain
$$
    I \supset J^{(1)} \supset J^{(2)} \supset \dots
$$
containes only a finite number of nonzero terms.
Appliing the formula~(\ref{rate_hilb}) to the ideals
$J^n$,
we obtain a finite sum
$$
    I(z) = z^{d(I)} (R(z) - N(I)(z))
            + z^{d(J^{(1)})} ((R(z) - N(J^{(1)})(z)))
            + z^{d(J^{(2)})} ((R(z) - N(J^{(2)})(z)))
	    + \dots
$$
Thus
$$
    I_i(z) = \sum_{j=1}^s a_{ij} \left( R(z) - I_j(z) \right),
$$
where $a_{ij} \in z {\bf Z} [z]$.
%= q_{ij} z^{d_{ij}}, q_{ij} \in {\bf Z}, d \ge d_{ij} \ge 1$.

Let $H = H(z)$ be a  column vector $[I_1(z), \dots, I_s(z)]^{t}$,
let $A$ be the matrix $(a_{ij}) \in z M_s({\bf Z}[z])$, and let
$e$ be a unit $s$--dimensional column vector.
We have
$$
    H = A \left( R(z)e - H \right),
$$
or
$$
   \left( A + E \right) H = R(z) Ae,
$$
where $E$ is the unit matrix.

The determinant $D(z) = \det \left( A + E \right) \in {\bf Z}[z]$
is a polynomial of degree at most $sd$. It is
invertible in ${\bf Q} \left[ [z] \right]$.
Then $\left( A + E \right) = D(z)^{-1} B$ with $B \in M_s({\bf Z}[z])$.
The elements of $B$ are minor determinants of $\left( A + E \right)$,
so, their degrees do not exceed $d(s-1)$.
We have
$$
    H = R(z) D(z)^{-1} C,
$$
where $C = BAe \in z {\bf Z}[z]^s$.
So, for every $1 \le i \le s$
$$
    I_i(z) = R(z) C_i(z) D(z)^{-1}.
$$

Let $I_s = \overline R$. Then
$$
   \overline R (z) = R(z) -1 = R(z) C_s(z) D(z)^{-1},
$$
therefore,
$$
     R(z) = \frac{\d D(z)}{\d D(z) - C_s(z)},
$$
so, $R(z)$ is a qoutient of two polynomials of degrees at most $sd$.
By the above,
$$
    I_i(z) = R(z) C_i(z) D(z)^{-1} = \frac{\d C_i(z)}{\d D(z) - C_s(z)},
$$
so, the same is true for the Hilbert series $I(z)$.

\begin{corollarium}
Suppose an algebra $R$ has a
Koszul filtration ${\bf F}$.
Then Hilbert series and  Poincare series of $R$
and of any ideal $I \in {\bf F}$
are rational functions.
\end{corollarium}

The statement about Poincare series follows from the equalities
$P_R(-z) R(z) = -1$ and
$P_I(-z) R(z) = I(z)$.

\section{Open questions}

\label{sectquest}

\begin{quest}
Does any commutative algebra with quadratic Groebner basis
of commutative relations ({\it $G$--quadratic algebra})
contain a Koszul flag?
\end{quest}

\begin{quest}
Are there PBW algebras which do not  have
Koszul filtrations?
\end{quest}

%Proposed answer is {\it Yes}.

\begin{quest}
Does any algebra having an infinite Koszul filtration has also a finite one?
\end{quest}

\begin{quest}
Suppose that an algebra $R$ has a Koszul filtration (or a Koszul flag).
Is the same true for its dual algebra $R^!$?
\end{quest}

\begin{quest}
Suppose that an algebra $R$ has a Koszul filtration or a Koszul flag
of right-sided ideals.
Is the same true for left-sided ideals?
\end{quest}

Notice that, at least, Koszul filtration of right-sided ideals
may be not of the same kind as the filtration of left-sided ideals.
For example, in the algebra $A = k \langle x,y,z,t |  zy-tz, zx \rangle $
{\it all} degree-one generated {\it left--}sided ideals
form a Koszul filtration
(so, $A$ is {\it universally Koszul}~\cite{conca} for left ideals), but
{\it right--}sided ideal $zA$ is not Koszul
module~\cite[Proposition 10]{pi2}.
However, $A$ is initially Koszul for
right-sided ideals.

\begin{quest}
Does any algebra with (infinite) rate filtration has rational
Hilbert series?
Does the same true for its Poincare series
$P_R(s,t) = \sum\limits_{i,j} \dim H_{ij}(R) s^i t^j$ ?
\end{quest}

\begin{quest}
It is shown by Backelin~\cite{brate} that, for any algebra $R$
of finite rate, its Veronese subalgebras (that is,
algebras of the type $R^{(d)} = \oplus_i R_{di}$) of sufficiently
high order $d$
are Koszul. Suppose that $R$ has a rate filtration of finite degree.
Do its Veronese subalgebras of high order have Koszul filtrations?
\end{quest}

\begin{quest}
Are universally Koszul
or ``universally rate'' algebras
always (projective) coherent?
%It is shown in~\cite{pi2} that algebras of restricted processing are
%always coherent. Is the same true for ?
\end{quest}


\begin{thebibliography}{Bib}




\bibitem[An1]{an1} D.~Anick,  {\it Generic algebras and CW--complexes},
                  Proc. of 1983 Conf. on algebra, topol. and K--theory
                  in honor of John Moore. Princeton Univ., 1988, p.~247--331
\bibitem[An2]{an2} D.~Anick,  {\it On the homology of associative algebras},
                  Trans. Amer. Math. Soc., {\bf 296} (1986), 2, p.~641--659
\bibitem[ABH]{abh} A. Aramova, S. Barcanescu, and J. Herzog,
		  {\it On the rate of relative Veronese submodules},
		  Rev. Roumaine Math. Pure Appl., {\bf 40} (1995),
		  3--4, p.~243--251
%\bibitem[AS]{AS}  M.~Artin, W. Shelter, {\it Graded algerbras of global
%                  dimension 3,} Adv. Math., {\bf 66} (1987), 171--216
%\bibitem[ATB]{ATB} M.~Artin, J.~Tate, M.~van~den~Bergh, {\it Some
%                   algebras related to authomorphisma of elliptic curves,}
%	           {\it The Grothendieck Festschrift,} v.1, 33--85, Birkhauser,
%                   Boston,1990
\bibitem[B]{burbaki} N. Bourbaki, {\it Alg\`ebre, ch.~X.
     Alg\`ebre homologique,} 1980
\bibitem[Ba1]{bac} Backelin~J., {\it A distributiveness property of augmented
  algebras, and some related homological results,} Ph. D. thesis, Stockholm
  (1982)
\bibitem[Ba2]{brate} Backelin~J., {\it On the rates of growth of homologies
		   of Veronese subrings}, Lecture Notes Math.,
                   {\bf 1183} (1986), 79--100
\bibitem[Bl]{i-kos} S. Blum, {\it Initially Koszul algebras},
                   Beitr\"age Algebra Geom., {\bf 41} (2000), 2, p.~455--467
\bibitem[BGS]{bgs} A.~Beilinson, V.~Ginzburg, W.~Soergel,
		   {Koszul duality patterns in representation theory,}
		   J. AMS, {\bf 9} (1996), 2, p.~473--527
\bibitem[Co1]{conca} A.~Conca, {\it Universally Koszul algebras}, Math. Ann.,
		   {\bf 317} (2000), 2, p.~329--346
\bibitem[Co2]{con2} A.~Conca, {\it Universally Koszul algebras
                   defined by monomials},
                   Rendiconti del seminario matematico
                   dell'Universita di Padova, {\bf 107} (2002), p.~1--5
\bibitem[CNR]{cnr} Conca A., de Negri E., and Rossi M. E.,
		   {\it On the rate of points in projective spaces},
		   Israel J. Math., {\bf 124} (2001), p.~253--265
\bibitem[CRV]{crv} Conca A., Rossi M. E., and Valla G.,
                   {\it Groebner flags and Gorenstein algebras,}
                   Compositio Math., {\bf 129} (2001), p.~95--121
\bibitem[CTV]{ctv} Conca A., Trung N. V., and Valla G.,
                   {\it Koszul property for
                   points in projective spaces,} Math. Scand.,
                   {\bf 89} (2001), 2, p.~201--216
\bibitem[D]{dav}   Davydov A. A., {\it Totally positive sequences
                   and R-matrix quadratic algebras,}
		   J. Math. Sci. (New York), {\bf 100} (2000), 1,
                   p.~1871--1876
\bibitem[F]{faith} K. Faith, {\it Algebra: rings, modules, and categories},
     v.~1
\bibitem[G]{gov}   Govorov V. E., {\it On graded algebras}, Math. Notes,
		   {\bf 12} (1972), 2, p.~197--204
\bibitem[HHR]{HHR} Herzog J., Hibi T., and Restuccia G.,
                   {\it Strongly Koszul algebras}, Math. Scand.,
		   {\bf 86} (2000), 2, p.~161--178
\bibitem[I]{iou}   Ioudu N. K.., {\it Solvability of the problem of
                   zero-divisor determining in a class of algebras},
		   Fundamentalnaja i prikladnaja matematika,
		   {\bf 1} (1995), 2, p.~541--544 [in Russian]
\bibitem[Pi1]{pion} Piontkovski D. I., {\it On the Hilbert series
                   of Koszul algebras,}
                   Functional Anal. and its Appl., {\bf 35} (2001), 2
\bibitem[Pi2]{pi2} Piontkovski D. I., {\it Noncommutative Groebner bases,
		   coherence of associative algebras, and devideness
                   in semigroups,} Fundamentalnaja i prikladnaja matematika,
		   {\bf 7} (2001), 2, p.~495--513 [in Russian]
\bibitem[Pi3]{pi3} Piontkovski D. I., {\it Hilbert series and relations
		   in algebras,} Izvestia: Mathematics, {\bf 6} (2000),
		   64, p.~1297--1311
\bibitem[PP]{pp}   Polishchuk A. and Positselski L.,
                   {\it Qudratic algebras, } preprint, 1996
\bibitem[Pol]{pol} Polishchuk A.,
                   {\it Noncommutative Proj and coherent algebras, }
                   preprint arXiv:math.RA/0212182, 2002
\bibitem[Pos]{pos} Positselski L. E., {\it Relation between the Hilbert series
                  of dual quadratic algebras does not imply Koszulity,}
                  Functional Anal. and its Appl., {\bf 29,} 3, 83--87 (1995)
\bibitem[Pr]{pri}   Priddy S. B., {\it Koszul resolutions,}
                   Trans. AMS, {\bf 152} (1970), 1, p.~39--60
\bibitem[Ro]{roos} Roos J.--E., {\it On the characterization of Koszul algebras.
   Four counter-examples.},
   Comptes Rendus Acad.\ Sci.\ Paris, S\'er. I,
   {\bf 321,} 1, 15--20 (1995)
\bibitem[U1]{uf1}  V.~A.~Ufnarovsky, {\it Algebras defined by two
		   quadratic relations}, Researchs on theory of
		   rings, algebras, and modules, Mathematical researchs,
		   Kishinev, 76 (1984), p.~148--172 [Russian]
\bibitem[U2]{ufn}  V.~A.~Ufnarovsky, {\it Combinatorial and asymptotical
                   methods in algebra,} Sovr. probl. mat., Fund. napr.,
	           {\bf 57} (1990), p.~5--177 [Russian]
	           Engl. transl.: Algebra VI, Encycl. Math. Sci., Springer,
                   Berlin 1995, p.~1--196
\bibitem[W]{wam}   M. Wambst, {\it Complexes de Koszul quantiques},
		   Ann. Inst. Fourier, Grenoble, {\bf 43} (1993), 4,
		   p.~1089--1156
\bibitem[Z]{zhang} J.~J.~Zhang,
                   {\it Non-Noetherian regular rings of dimension 2,}
	           Proc. AMS, {\bf 126} (1998), 6, 1645--1653



% ========= skipped ==========

%\bibitem[AS]{AS}  M.~Artin, W. Shelter, {\it Graded algerbras of global
%                  dimension 3,} Adv. Math., {\bf 66} (1987), 171--216
%\bibitem[ATB]{ATB} M.~Artin, J.~Tate, M.~van~den~Bergh, {\it Some
%                   algebras related to authomorphisma of elliptic curves,}
%	           {\it The Grothendieck Festschrift,} v.1, 33--85, Birkhauser,
%                   Boston,1990
\end{thebibliography}
\end{document}